%% file: vapnik.tex
\documentclass[]{amsart}
\usepackage{fancyhdr}
\usepackage{graphicx}
\usepackage{amssymb}
\usepackage{dsfont}
\usepackage{eucal}
\usepackage{mathrsfs}
\newcommand{\ds}{\displaystyle}

\renewcommand{\mathbb}{\mathds}

\newcommand{\PP}{\mathbb{P}}

\newcommand{\RR}{\mathbb{R}}
\newcommand{\NN}{\mathbb{N}}

\newcommand{\C}[1]{\mathcal{#1}}

\newcommand{\rr}{\overline{r}}

\newcommand{\B}[1]{\mathds{#1}}

\DeclareMathOperator*{\ess}{ess}

\DeclareMathOperator{\Var}{\mathbb{V}ar}

\numberwithin{equation}{section}
\newtheorem{thm}{Theorem}[section]

\newtheorem{proposition}[thm]{Proposition}
\newtheorem{lemma}[thm]{Lemma}
\newtheorem{cor}[thm]{Corollary}

\theoremstyle{definition}
\newtheorem{dfn}{Definition}[section]
\theoremstyle{remark}
\newtheorem{rmk}{Remark}[section]

\fancypagestyle{plain}{
\fancyhf{}
\fancyhead[LO]{\tiny\rightmark}
\fancyhead[RE]{\tiny\leftmark}
\fancyhead[LE,RO]{\thepage}
\fancyfoot[RO,LE]{\tiny\sc Improved Vapnik Cervonenkis bounds}
\fancyfoot[LO,RE]{\tiny\sc Olivier Catoni -- \today}
\renewcommand{\footruleskip}{0pt}
\renewcommand{\footrulewidth}{0.4pt}
}

\allowdisplaybreaks
\begin{document}
\pagestyle{fancy}
\fancyhf{}
\fancyhead[LO]{\tiny\rightmark}
\fancyhead[RE]{\tiny\leftmark}
\fancyhead[LE,RO]{\thepage}
\fancyfoot[RO,LE]{\tiny\sc Improved Vapnik Cervonenkis bounds}
\fancyfoot[LO,RE]{\tiny\sc Olivier Catoni -- \today}
\renewcommand{\footruleskip}{0pt}
\renewcommand{\footrulewidth}{0.4pt}
%\pagestyle{headings}
%\frontmatter
\title{Improved Vapnik Cervonenkis bounds}
\author{Olivier Catoni}
\address{CNRS --   
Laboratoire de Probabilit\'es et Mod\`eles Al\'eatoires,
Universit\'e Paris 6 (site Chevaleret), 4 place Jussieu -- Case 188, 
75 252 Paris Cedex 05.}

\begin{abstract}
We give a new proof of VC bounds where we avoid the use of symmetrization
and use a shadow sample of arbitrary size. We also improve 
on the variance term. This results in better
constants, as shown on numerical examples. Moreover our bounds
still hold for non identically distributed independent random variables.\\[12pt]
\end{abstract}
\maketitle
%\tableofcontents
%\mainmatter
%{\fontencoding{U}\fontfamily{ygoth}\selectfont}
%\include{introBis}
\input{text}

\input{biblio}
\end{document}

%% file: text.tex
%$Id: text.tex,v 1.2 2004/10/10 09:59:34 catoni Exp catoni $
\section{Description of the problem}
Let $(\C{X}, \C{B})$ be some measurable space and $\C{Y}$ some
finite set. Let $(\Theta, \C{T})$ be a measurable parameter space
and $\bigl\{ f_{\theta} : \C{X} \rightarrow \C{Y}, 
\theta \in \Theta \}$ be a family of decision functions. Assume that
$$
(\theta, x) \mapsto f_{\theta}(x) : ( \Theta \times \C{X}, 
\C{T} \otimes \C{B} ) \rightarrow \C{Y}
$$
is measurable.
Let 
$$
P_i \in \C{M}_+^1\bigl(\C{X} \times \C{Y}, \C{B} \otimes \{ 0,1\}^{\C{Y}}
\bigr), \quad i=1, \dots, N,
$$ be some probability distributions on $\C{X} \times \C{Y}$
--- where $\{0,1\}^{\C{Y}}$ is the discrete sigma algebra of 
all the subsets of $\C{Y}$.
Let $(X_i, Y_i)_{i=1}^N$ be the canonical process on $(\C{X} \times \C{Y})^N$
--- i.e. the coordinate process $(X_i, Y_i)(\omega) = \omega_i$, $\omega
\in \bigl(\C{X} \times \C{Y})^N$. Let
$$
r(\theta) = \frac{1}{N} \sum_{i=1}^N \B{1}\bigl[ f_{\theta}(X_i)
\neq Y_i \bigr].
$$
We are interested in bounding with $\bigotimes_{i=1}^N P_i$ probability
at least $1-\epsilon$ and for any $\theta \in \Theta$
the quantity $R(\theta) - r(\theta)$. This question has an interest
both in statistical learning theory and in empirical process 
theory.

In the case when $ \lvert \C{Y} \rvert = 2$, introducing the notation
$$
\C{N}(X_1^{2N}) = \bigl\lvert \bigl\{ \bigl[ f_{\theta}(X_i) \bigr]_{i=1}^{2N};
\theta \in \Theta \bigr\} \bigr\rvert,
$$
where $\lvert A \rvert$ is the number of elements of the set $A$,
Vapnik proved in \cite[page 138]{Vapnik} that
\begin{thm}
\label{th1}
For any probability distribution $P \in \C{M}_+^1
\bigl( \C{X} \times \C{Y} \bigr)$, with $P^{\otimes N}$ probability at least $1 - \epsilon$, 
for any $\theta \in \Theta$, 
$$
R(\theta) \leq r(\theta) + \frac{2 d'}{N} \biggl( 1 + 
\sqrt{\ds 1 + \frac{Nr(\theta)}{d'}} \biggr),
$$
where
$$
d' = \log \Bigl\{  P^{\otimes 2N} \bigl[ \C{N}(X_1^{2N} ) \bigr] 
\Bigr\} + \log \bigl( 4 
\epsilon^{-1} \bigr).
$$
\end{thm}
It is also well known since the works of Vapnik and Cervonenkis
that, in the case when $\C{Y} = \{0,1\}$, 
$$
\log \bigl[ \C{N}(X_1^{2N}) \bigr] \leq h \log \left( \frac{eN}{h} 
\right),
$$
where
$$
h = \max \bigl\{ \lvert A \rvert; \C{N}\bigl[(X_i)_{i \in A} \bigr] = 
2^{\lvert A \rvert} \bigr\}.
$$
Therefore when the VC dimension of $\{ f_{\theta}; \theta \in 
\Theta \}$ is not greater than $h$, that is when by definition
$$
\max \Bigl\{ \lvert A \rvert; A \subset \C{X}, 
\bigl\lvert \bigl\{ \bigl( f_{\theta}(x) \bigr)_{x \in A}; 
\theta \in \Theta \bigr\} \bigr\rvert = 2^{\lvert A \rvert}
\Bigr\} \leq h,
$$
we have the following
\begin{cor}
\label{cor1}
When the VC dimension of $\bigl\{ f_{\theta}; \theta \in \Theta \}$
is not greater than $h$, with $P^{\otimes N}$ probability at least
$1 - \epsilon$, for any $\theta \in \Theta$,
$$
R(\theta) \leq r(\theta) + \frac{2 d'}{N} \biggl(
1 + \sqrt{\ds 1 + \frac{N r(\theta)}{d'}} \biggr),
$$
where
$$
d' = h \log \biggl(\frac{2eN}{h} \biggr) + \log \bigl( 4 \epsilon^{-1} \bigr).
$$
\end{cor}
The aim of this paper is to improve theorem \ref{th1} and its corollary,
using PAC-Bayesian inequalities with data dependent priors.

We have already proved in \cite{CatPAC} that with $P^{\otimes N}$ probability 
at least $1 - \epsilon$, 
\begin{equation}
\label{previousBound}
R(\theta) \leq r(\theta) + \frac{\zeta d}{N}
\left( 1 + \sqrt{\ds 1 + \frac{4 N r(\theta)}{\zeta d}}\, \right),
\end{equation}
where 
$$
d = P^{\otimes N}_{X_{N+1}^{2N}} \Bigl\{ 
\log \bigl[ \C{N}(X_1^{2N})\bigr] \Bigr\} + \log \biggl(\frac{\log(2\zeta N)}{\epsilon \log(
\zeta)} \biggr),
$$
which brings an improvement when $r(\theta) \leq \frac{d}{N}$ and 
$d$ is large.

Here we are going to generalize this theorem to arbitrary shadow sample
sizes and non identically distributed independent random variables.
We will also improve on the variance term in \eqref{previousBound} 
and get rid of the (unwanted !) parameter $\zeta$.

Moreover, we will derive VC bounds in the transductive setting in which the
shadow sample error rate is bounded in terms of the empirical
error rate (in this setting the shadow sample would more 
appropriately be described as a test set).

We will start with the transductive setting, since it has an interest
of its own and will in the same time serve as a 
technical step towards more classical results.

\section{The transductive setting}

We will consider a shadow sample of size
$kN$ where $k$ is some integer.\\
Let $(X_i, Y_i)_{i=1}^{(k+1)N}$
be the canonical process on $\bigl( \C{X} \times \C{Y} \bigr)^{(k+1)N}$.\\
We assume that we observe the first sample $(X_i, Y_i)_{i=1}^N$, 
that we may also observe the rest of the design 
$X_{N+1}^{(k+1)N}$, (this is a short notation
for $(X_i)_{i=N+1}^{(k+1)N}$), but that we
do not observe $Y_{N+1}^{(k+1)N}$.\\
Let $r_1(\theta)$ and $r_2(\theta)$ be the empirical
error rates of the decision function $f_{\theta}$ on 
the training and test sets:
\begin{align*}
r_1 (\theta) & = \frac{1}{N} \sum_{i=1}^N \B{1}\bigl[ Y_i \neq f_{\theta}(X_i) \bigr],\\
r_2 (\theta) & = \frac{1}{kN} \sum_{i=N+1}^{(k+1)N} 
\B{1} \bigl[ Y_i \neq f_{\theta}(X_i) \bigr].
\end{align*}
Let $\PP \in \C{M}_+^1\Bigl[ \bigl(\C{X} \times \C{Y}\bigr)^{(k+1)N} 
\bigr]$ be some {\em partially exchangeable} probability distribution on 
$\bigl( \C{X} \times \C{Y} \bigr)^{(k+1)N}$. What we mean by {\em partially 
exchangeable} will be precisely defined in the following. An important case
is when $\B{P} = \left( \bigotimes_{i=1}^N P_i \right)^{\otimes (k+1)}$, 
meaning that we have $(k+1)$ independent samples, each being distributed
according to the same product of non identical probability distributions.
Let as in the introduction
$$
\C{N}(X_1^{(k+1)N}) = \Bigl\lvert \Bigl\{ \bigl[f_{\theta}(X_i) \bigr]_{i=1}^{(k+1)N} :
\theta \in \Theta \Bigr\} \Bigr\rvert
$$
be the number of distinct decision rules induced by the model on the design
$(X_i)_{i=1}^{(k+1)N}$.
We will prove
\begin{thm}
\label{th2.1}
With $\B{P}$ probability at least $1 - \epsilon$, 
for any $\theta \in \Theta$, 
$$
r_2(\theta) \leq r_1(\theta) + \frac{d}{N} + 
\sqrt{\ds \frac{2d(1+\frac{1}{k}) r_1(\theta)}{N} + \frac{d^2}{N^2}},
$$
where $d = \log \bigl[ \C{N}(X_1^{(k+1)N}) \bigr] + \log(\epsilon^{-1})$.
\end{thm}
Let us remind that when $\lvert \C{Y} \rvert = 2$ and the VC dimension of 
$\{ f_{\theta}; \theta \in \Theta \}$ is not greater than $h$, 
$$
d \leq h \log \biggl( \frac{e(k+1) N}{h} \biggr) + \log(\epsilon^{-1}).
$$

Let us take some numerical example : when $N = 1000$, $h = 10$,
$\epsilon = 0.01$ and $r_1(\theta) = 0.2$, we get
$r_2(\theta) \leq 0.4872$ using $k = 4$ (whereas for $k=1$
we get only $r_2(\theta) \leq 0.5098$, showing that increasing
the shadow sample size is useful to get a bound less than $0.5$) 

\newcommand{\Pk}{\PP}
Let us start the proof of theorem \ref{th2.1} with some notations and 
a few lemmas.
Let $\chi_i = \B{1} \bigl[ Y_i \neq f_{\theta}(X_i) \bigr] 
\in \{0,1\}$. For any random variable $h : \Omega = \bigl( \C{X} \times \C{Y} \bigr)^{(k+1)N} 
\rightarrow \RR$ ( we work on the
canonical space), let the transformed random variable $\tau_i (h)$ be 
defined as 
$$
\tau_i (h) = \frac{1}{k+1} \sum_{j=0}^k
h \circ \tau_i^j,
$$
where $\tau_i^j : \Omega \rightarrow \Omega$ is defined by 
$$
\bigl[\tau_i^j(\omega)\bigr]_{\ell} = 
\begin{cases}
\omega_{i+mn}, & \ell =i+ [(m + j) \bmod (k+1)] N, \quad m=0, \dots, k;\\
\omega_{\ell},& \ell \not\in \{i+mN : m=0, \dots, k \}.
\end{cases}
$$
In other words, $\tau_i^j$ performs a circular permutation of the 
subset of indices \linebreak $\{ i + m N : m=0, \dots, k\}$.
Notice also that $\tau_i$ may be viewed as a regular conditional
probability measure.
\begin{dfn}
The joint distribution $\PP$ is said to be partially exchangeable
when for any $i=1, \dots, N$, any $j = 0, \dots, k$, 
$\PP \circ (\tau_i^j)^{-1} = \PP$.
\end{dfn}
Equivalently, this means that for any bounded random variable $h$, 
$$
\PP(h) = \PP( h \circ \tau_i^1), \quad i=1, \dots, N,
$$
(since $\tau_i^j$ is the $j$th iterate of $\tau_i^1$).
As a result, any partially exchangeable distribution $\PP$
is such that for any bounded random variable
$$
\Pk(h) = \Pk \left\{ \left[ \bigcirc_{i=1}^N \tau_i \right] (h) \right\},
$$
where we have used the notation $\bigcirc_{i=1}^N \tau_i = \tau_1 \circ \tau_2
\circ \dots \circ \tau_N$.

In the same way
\begin{dfn} A random variable $h : \bigl(\C{X} \times \C{Y} \bigr)^{(k+1)N} 
\rightarrow \RR$ is said to be partially exchangeable when for any 
$i = 1, \dots, N$, $h \circ \tau_i^1 = h$.
\end{dfn}

\begin{lemma}
\label{lemma2.1}
For any $\theta \in \Theta$, 
any $\omega \in \bigl( \C{X} \times \C{Y} \bigr)^{(k+1)N}$, 
any positive partially exchangeable random variable $\lambda$,
any partially exchangeable random variable $\eta$,
$$
\Bigl( \bigcirc_{i=1}^N \tau_i \Bigr) \Bigl\{ \exp \Bigl[ \lambda \bigl[ r_2(\theta)
- r_1(\theta) \bigr] - \eta \Bigr] \Bigr\} (\omega) 
\leq \exp \Bigl[ \tfrac{\lambda^2}{2 N} \bigl[ \tfrac{1}{k} r_1(\theta) 
+ r_2(\theta)\bigr] - \eta \Bigr] (\omega).
$$
\end{lemma} 
\begin{proof}
\begin{multline*}
\Bigl( \bigcirc_{i=1}^N \tau_i \Bigr) \Bigl\{ \exp \Bigl[ \lambda \bigl[ 
r_2(\theta) - r_1(\theta) \bigr] - \eta \Bigr] \Bigr\}
\\ = \exp ( - \eta) \prod_{i=1}^N \tau_i \biggl\{ \exp 
\biggl( \frac{\lambda}{kN} \sum_{j=1}^k \chi_{i+jN} - \frac{\lambda}{N} 
\chi_i \biggr) \biggr\}\\
= \exp(-\eta) \prod_{i=1}^N \exp \biggl( \frac{\lambda}{kN}
\sum_{j=0}^k \chi_{i+jN} \biggr) 
\prod_{i=1}^N \tau_i \biggl\{ \exp \biggl( - \frac{(k+1)\lambda}{kN} \chi_i \biggr)
\biggr\}.
\end{multline*}
Let $p_i = \frac{1}{k+1} 
\sum_{j=0}^k \chi_{i+jN}$. Let $\chi$ be the identity (seen as the 
canonical process) on $\{0,1\}$ and 
$B_p$ be the Bernoulli distribution on $\{0, 1\}$ with parameter $p$,
namely let $B_p(1) = 1 - B_p(0) = p$. It is easily seen that
$$
\log \biggl\{ \tau_i \biggl[ \exp \biggl( 
- \frac{(k+1)\lambda}{kN} \chi_i \biggr) \biggr] \biggr\} 
= \log \biggl\{ B_{p_i} \biggl[ \exp \biggl(
- \frac{(k+1)\lambda}{kN} \chi \biggr) \biggr] \biggr\}.
$$
Moreover this last quantity can be bounded in the following way.
$$
\log \Bigl\{ B_p \bigl[ \exp (- \alpha \chi) \bigr] \Bigr\}
= -\alpha B_p(\chi) + \int_0^{\alpha} (1 - \beta) \Var_{B_{f(\beta)}}(
\chi) d \beta.
$$
This is the Taylor expansion of order two of $\alpha \mapsto \log \Bigl\{ 
B_p \bigl[ \exp ( -\alpha \chi ) \bigr] \Bigr\}$, 
where 
$$
f(\beta) = \frac{B_p\bigl[ \chi \exp ( - \beta \chi)\bigr]}{
B_p\bigl[\exp(- \beta \chi)\bigr]} = 
\frac{p \exp(-\beta)}{(1-p)+ p \exp(-\beta)} \leq p.
$$
Thus
$$
\Var_{B_{f(\beta)}}(\chi) = f(\beta)\bigl[ 1 - f(\beta) \bigr] 
\leq (p \wedge \tfrac{1}{2})\bigl[ 1 - (p \wedge \tfrac{1}{2}) \bigr]
\leq (1 - \tfrac{1}{k+1}) p = \tfrac{k}{k+1} p,
$$
for any $p \in \bigl( \frac{1}{k+1} \NN \bigr) \cap [0,1]$.
Hence
$$
\log \Bigl\{ B_p \bigl[ \exp ( - \alpha \chi ) \bigr] \Bigr\} 
\leq - \alpha p + \frac{k \alpha^2}{2(k+1)} p,
$$
and 
$$
\tau_i \biggl[ \exp \biggl( - \frac{(k+1)\lambda}{kN} \chi_i \biggr) \biggr]
\leq \exp \biggl( - \frac{(k+1)\lambda}{kN} p_i + \frac{(k+1)\lambda^2}{2kN^2} 
p_i \biggr).
$$
Therefore
\begin{multline*}
\Bigl( \bigcirc_{i=1}^N \tau_i \Bigr)  \biggl\{ \exp \Bigl[ \lambda 
\bigl[ r_2(\theta) - r_1(\theta) \bigr] - \eta \Bigr] \biggr\} 
\leq \exp ( - \eta) \exp \Biggl( \frac{(k+1) \lambda^2}{2kN^2} 
\sum_{i=1}^N p_i \Biggr) \\
= \exp \Biggl( \frac{\lambda^2}{2kN^2} \sum_{i=1}^{(k+1)N} 
\chi_i - \eta \Biggr) = \exp \biggl\{ \frac{\lambda^2}{2N} 
\bigl[ \tfrac{1}{k}r_1(\theta) + r_2(\theta) \bigr] - \eta \biggr\}.
\end{multline*}
\end{proof}

\begin{lemma}
For any $\theta \in \Theta$, for any positive partially 
exchangeable random variable $\lambda$, 
for any partially exchangeable random variable
$\eta$, 
$$
\Pk \Bigl\{ \exp \Bigl[ \lambda \bigl[ r_2(\theta) - r_1(\theta) \bigr] 
- \eta \Bigr] \Bigr\} \leq \Pk \Bigl\{ 
\exp \Bigl[ \tfrac{\lambda^2}{2N} \bigl[ \tfrac{1}{k}r_1(\theta) + r_2(\theta)
\bigr] - \eta \Bigr] \Bigr\}.
$$
\end{lemma}
\begin{rmk}
Let us notice that we do not need integrability conditions, and that
the previous inequality between expectations of positive random variables
holds in $\RR_+ \cup \{+ \infty\}$, meaning that both members may be equal
to $+ \infty$.
\end{rmk}
\begin{rmk}
We can take $\eta = \log(\epsilon^{-1}) + \frac{\lambda^2}{2N} \bigl[
\frac{1}{k}r_1(\theta) + r_2(\theta) \bigr]$ to get 
$$
\Pk \Bigl\{ \exp \Bigl[ \lambda \bigl[ r_2(\theta) - r_1(\theta) \bigr] 
 - \tfrac{\lambda^2}{2N} \bigl[ \tfrac{1}{k} r_1(\theta) + r_2(\theta)
 \bigr] + \log(\epsilon) \Bigr] \Bigr\} \leq \epsilon.
$$
\end{rmk}
\begin{proof}
According to the previous lemma,
\begin{multline*}
\Pk \Bigl\{ \exp \Bigl[ \lambda \bigl[ r_2(\theta) - r_1(\theta) 
\bigr] - \eta \Bigr] \Bigr\} \\ = 
\Pk \biggl\{ \Bigl( \bigcirc_{i=1}^N \tau_i \Bigr) \Bigl\{ \exp \Bigl[ 
\lambda \bigl[ r_2(\theta) - r_1( \theta) \bigr] - \eta \Bigr] \Bigr\}
\biggr\} \\
\leq \Pk \biggl\{ \exp \biggl[ \frac{\lambda^2}{2N} \bigl[ 
\tfrac{1}{k} r_1(\theta) + r_2(\theta) \bigr] - \eta \biggr] \biggr\}.
\end{multline*}
\end{proof}

Let us now consider some {\em partially exchangeable prior distribution} 
$\pi \in \C{M}_+^1(\Theta)$:
\begin{dfn}
A regular conditional probability distribution\\ $\pi : \bigl(\C{X} \times
\C{Y} \bigr)^{(k+1)N} \rightarrow \C{M}_+^1(\Theta, \C{T})$ is said to 
be partially exchangeable when for any $i=1, \dots, N$, any $\omega \in 
\bigl(\C{X} \times \C{Y} \bigr)^{(k+1)N}$, 
$\pi\bigl[\tau_i^1(\omega)\bigr] = \pi(\omega)$,
this being an equality between probability measures in $\C{M}_+^1(\Theta,\C{T})$.
\end{dfn}
In the following, $\lambda$ and $\eta$ will be random variables depending
on the parameter $\theta$. We will say that a real random variable
$h : \bigl(\C{X} \times \C{Y} \bigr)^{(k+1)N} \times \Theta \rightarrow \RR$ is 
partially exchangeable when $h(\omega, \theta) = 
h\bigl[\tau_i^1(\omega), \theta \bigr]$, $i=1, \dots, N$, $\theta \in \Theta$,
$\omega \in \bigl( \C{X} \times \C{Y} \bigr)^{(k+1)N}$.
\begin{lemma} For any partially exchangeable prior distribution $\pi$,
any positive partially exchangeable random variable $\lambda : 
\bigl( \C{X} \times \C{Y} \bigr)^{(k+1)N} \times \Theta \rightarrow \RR$, 
and any partially exchangeable random threshold function 
$\eta : \bigl( \C{X} \times \C{Y} \bigr)^{(k+1)N} \times \Theta \rightarrow \RR$,
\begin{multline*}
\Pk \biggl\{ \pi \Bigl\{ \exp \Bigl[ \lambda(\theta) \bigl[
r_2(\theta) - r_1(\theta) \bigr] - \eta(\theta) \Bigr] \Bigr\} \biggr\} 
\\ \leq \Pk \Biggl\{ \pi \biggl\{ \exp \biggl[ 
\frac{\lambda(\theta)^2}{2N} \bigl[ \tfrac{1}{k}r_1(\theta) + r_2(\theta) \bigr]
- \eta(\theta)
\biggr] \biggr\} \Biggr\}.
\end{multline*}
\end{lemma}
\begin{proof} It is a consequence of lemma \ref{lemma2.1} 
and of the following identities:
\begin{multline*}
\Pk \biggl\{ \pi \Bigl\{ 
\exp \Bigl[ \lambda(\theta) \bigl[ r_2(\theta) - r_1(\theta) 
\bigr] - \eta(\theta) \Bigr] \Bigr\} \biggr\} 
\\ = \Pk \biggl\{ \Bigl( \bigcirc_{i=1}^N 
\tau_i \Bigr) \biggl( \pi \Bigl\{ 
\exp \Bigl[ \lambda(\theta) \bigl[ r_2(\theta) - r_1(\theta) 
\bigr] - \eta(\theta) \Bigr] \Bigr\} \biggr) \biggr\}\\
= \Pk \biggl\{ 
\pi \Bigl\{ \bigl( \bigcirc_{i=1}^N \tau_i \bigr) 
\exp \Bigl[ \lambda(\theta) \bigl[ r_2(\theta) - r_1(\theta) 
\bigr] - \eta(\theta) \Bigr] \Bigr\} \biggr\}.
\end{multline*}
Indeed for any positive random variable $h : \bigl( \C{X} \times \C{Y} 
\bigr)^{(k+1) N} \times \Theta \rightarrow \RR$, 
$$
\pi(h) \circ \tau_i^j = 
(\pi \circ \tau_i^j) (h \circ \tau_i^j) 
= \pi ( h \circ \tau_i^j ).
$$
Thus 
$$
\tau_i \bigl[ \pi(h) \bigr] = \frac{1}{k+1} \sum_{j=0}^k \pi \bigl(
h \circ \tau_i^j \bigr)  
= \pi \biggl( \frac{1}{k+1} \sum_{j=0}^k h \circ \tau_i^j 
\biggr) = \pi \bigl( \tau_i h \bigr).
$$
\end{proof}
As a consequence, we get the following learning theorem:
\begin{thm}
\label{th2.4}
For any partially exchangeable prior distribution $\pi$, 
any positive partially exchangeable random variable $\lambda$, 
with $\Pk$ probability at least $1 - \epsilon$, for any 
$\rho \in \C{M}_+^1(\Theta)$, 
$$
\rho \bigl[ \lambda (\theta) r_2(\theta) \bigr] 
- \rho \bigl[ \lambda( \theta) r_1(\theta) \bigr] 
\leq \rho \biggl\{ \frac{\lambda(\theta)^2}{2N} \Bigl[ 
\tfrac{1}{k} \bigl[ r_1(\theta) +  
r_2(\theta) \Bigr] \biggr\}  
+ \C{K}(\rho, \pi) + \log(\epsilon^{-1}).
$$
\end{thm}
\begin{proof}
Take $\eta (\theta) = \frac{\lambda(\theta)^2}{2N} 
\bigl[ \tfrac{1}{k} r_1(\theta) + r_2(\theta) 
\bigr] + \log(\epsilon^{-1})$ and notice that it 
is indeed a partially exchangeable threshold function. 

\pagebreak
Thus
\begin{multline*}
\Pk \biggl\{ \sup_{\rho \in \C{M}_+^1(\Theta)} 
\rho \bigl[ \lambda(\theta) r_2(\theta) \bigr] - \rho 
\bigl[\lambda(\theta) r_1(\theta) \bigr] \\ - \rho
\biggl[ \frac{\lambda(\theta)^2}{2N} \Bigl\{ \tfrac{1}{k} 
\rho\bigl[ r_1(\theta) \bigr] + \rho \bigl[ r_2(\theta) \bigr] 
\Bigr\} \biggr] - \C{K}(\rho,\pi) 
+ \log(\epsilon) \leq 0 \biggr\} 
\\ = 
\Pk \biggl\{ \log \Bigl\{ \pi \Bigl[ \exp \bigl\{ 
\lambda(\theta) \bigl[ r_2(\theta) - r_1(\theta) \bigr] \\ 
\shoveright{- \tfrac{\lambda(\theta)^2}{2N}
\bigl[ \tfrac{1}{k}r_1(\theta) + r_2(\theta) \bigr] + \log(\epsilon) \bigr\}
\Bigr]
\Bigr\} \leq 0 \biggr\}} \\
\leq \Pk \biggl\{ \pi \biggl[ 
\exp \Bigl\{ \lambda(\theta) \bigl[ r_2(\theta) - r_1(\theta) \bigr] 
\\ - \tfrac{\lambda(\theta)^2}{2N} \bigl[ \tfrac{1}{k} r_1(\theta) 
+ r_2(\theta) \bigr] + \log(\epsilon) \Bigr\} \biggr] \biggr\}
\leq \epsilon.
\end{multline*}
We have used the identity $\log \Bigl\{ \pi \bigl[
\exp ( h ) \bigr] \Bigr\} = \sup_{\rho \in \C{M}_+^1(\Theta)}
\rho(h) - \C{K}(\rho, \pi)$. See for instance 
\cite[pages 159-160]{Cat7} or \cite[lemma 4.2]{CatPAC} for a proof.
\end{proof}

Let us consider the map $\Psi : \Theta \rightarrow \C{Y}^{(k+1)N}$
which restricts each classification rule
to the design: 
$\Psi(\theta) = \bigl[ f_{\theta}(X_i) \bigr]_{i=1}^{(k+1)N}$.
Let $\Theta / \Psi$ be the set of components of $\Theta$
for the equivalence relation 
$\{ (\theta_1, \theta_2) \in \Theta^2; \Psi(\theta_1) = \Psi(\theta_2)
\}$. Let
$ c : \{0,1\}^{\Theta} \rightarrow \Theta$ be such that $c(\theta') \in 
\theta'$ for each $\theta' \subset \Theta$ (the function $c$ chooses
some element from any subset of $\Theta$).
Let $\Theta' = c \bigl(\Theta/\Psi \bigr)$. Let us note that
$\Psi$ and therefore $\Theta / \Psi$ and $\Theta'$ are exchangeable
random objects. Let
$$
\pi = \frac{1}{\lvert \Theta' \rvert} \sum_{\theta \in \Theta'} 
\delta_{\theta}
$$
be the uniform distribution on the finite subset $\Theta'$ of $\Theta$.

Applying theorem \ref{th2.4} to $\pi$,
and $\rho = \delta_{\theta}$,
we get that for any positive partially exchangeable 
random variable
$\lambda$, 
with $\PP$ probability at least $1 - \epsilon$, 
for any $\theta' \in \Theta'$, 
$$
\lambda(\theta') r_2(\theta') - \lambda(\theta') r_1(\theta') \leq 
\frac{\lambda(\theta')^2}{2N} 
\Bigl[ \tfrac{1}{k} r_1(\theta') + r_2(\theta') \Bigr] 
+ \log \bigl\lvert \Theta' \bigr\rvert
+ \log(\epsilon^{-1}).
$$
Let us choose 
$$
\lambda(\theta) = \left( \frac{2 N \log \Bigl(
\frac{\lvert \Theta' \rvert}{\epsilon} \Bigr)}{\frac{1}{k}r_1(\theta) + 
r_2(\theta)} \right)^{1/2},
$$
with the convention that when $\frac{1}{k}r_1(\theta) + r_2(\theta) = 0$,
then $\lambda r_2(\theta) = \lambda r_1(\theta) = 0$.
This is legitimate, since $\lvert \Theta' \rvert$ 
and $\frac{1}{k}r_1(\theta') 
+ r_2(\theta')$ are exchangeable random variables,
and since when $\frac{1}{k}r_1(\theta) + r_2(\theta) = 0$, then 
$r_1(\theta) = r_2(\theta) = 0$.

Thus, with $\PP$ probability at least $1 - \epsilon$, 
for any $\theta' \in \Theta'$, 
$$
r_2(\theta') - r_1(\theta') 
\leq \left( \frac{2 \log \bigl( \tfrac{\lvert \Theta' \rvert}{\epsilon} \bigr)
\bigl[ \frac{1}{k}r_1(\theta') + r_2(\theta') \bigr]}{N} \right)^{1/2}.
$$
Now we can remark that for each $\theta \in \Theta$, 
$\theta' = c \bigl[ \Psi(\theta) \bigr]$ is such that
$f_{\theta'}(X_i) = f_{\theta}(X_i)$, for $i = 1, \dots, (k+1)N$. 
Therefore $r_1(\theta) = r_1(\theta')$ and $r_2(\theta)
= r_2(\theta')$. 

Thus with $\PP$ probability at least $1 - \epsilon$,
for any $\theta \in \Theta$, 
\begin{equation}
\label{eq2.2}
r_2(\theta) - r_1(\theta) 
\leq \left( \frac{2 \log \bigl( \tfrac{\lvert \Theta' \rvert}{\epsilon} \bigr)
\bigl[ \frac{1}{k}r_1(\theta) + r_2(\theta) \bigr]}{N} \right)^{1/2}.
\end{equation}

Putting for short $d = \log \bigl( \frac{\lvert \Theta' \rvert}{\epsilon} \bigr)$
and solving inequality \eqref{eq2.2}
with respect to $r_2(\theta)$ proves theorem \ref{th2.1}.

Note that we have in fact proved a more general version of
theorem \ref{th2.1}, where $d$ can be taken to 
be $d = - \log \bigl[ \Bar{\pi}(\theta) \epsilon \bigr]$, where 
$$
\Bar{\pi}(\theta) = \sup \{ \pi(\theta') : \theta' \in \Theta, 
\Psi(\theta') = \Psi(\theta) \},
$$ for any choice of partially
exchangeable prior probability distribution $\pi$.

\section{Improvement of the variance term}

We will first improve the variance term in lemma \ref{lemma2.1}
when $k=1$, and $\PP$ is fully exchangeable. We will deal
afterwards with the general case.

\begin{thm}
\label{th3.3}
For any exchangeable probability distribution 
$\PP$, with $\PP$ probability $1 - \epsilon$, for any $\theta \in \Theta$, 
$$
r_2(\theta) \leq r_1(\theta) + \frac{d}{N}\bigl[ 1 - 2 r_1(\theta) \bigr]
+ \sqrt{\frac{4 d}{N} \bigl[1 - r_1(\theta) \bigr] r_1(\theta) + 
\frac{d^2}{N^2}\bigl[1 - 2 r_1(\theta) \bigr]^2 },
$$
where $d = \inf \Bigl\{ 
-\log\bigl[\pi(\theta') \epsilon \bigr]
: \theta' \in \Theta, \Psi(\theta') = \Psi(\theta) \Bigr\}$.
\end{thm}

Let us pursue our numerical example : assuming that $\lvert \C{Y} 
\rvert = 2$, $N = 1000$, $h = 10$,
$\epsilon = 0.01$
and $r_1(\theta) = 0.2$, we get that $r_2(\theta) \leq 0.453$.

\begin{proof}
Proving theorem \ref{th3.3} will require some lemmas.

Let
$$
\tau(h)(\omega) = \frac{1}{(2N)!} \sum_{\sigma \in \mathfrak{S}_{2N}}
h(\omega \circ \sigma), \quad \omega \in \Omega,
$$
where $\mathfrak{S}_{2N}$ is the set of permutations
of $\{1, \dots, 2N \}$ and where $(\omega \circ \sigma)_i = \omega_{\sigma(i)}$.
For any $\omega \in \Omega$, any $\sigma \in \mathfrak{S}_N$, 
let $\omega_{2,\sigma}$ be defined as 
$$
(\omega_{2,\sigma})_{i} = 
\begin{cases}
\omega_i, & 1 \leq i \leq N,\\
\omega_{\sigma(i-N)}, & N < i \leq 2N.
\end{cases}
$$
Let 
$$
\tau'(h)(\omega) = \frac{1}{N!} 
\sum_{\sigma \in \mathfrak{S}_N} h(\omega_{2,\sigma}).
$$
Let us remark that $\tau = \tau \circ \tau'$, and that
$\tau'\bigl[r_k(\theta)\bigr] = r_k(\theta)$, $k=1,2$.

Moreover, we know from the previous section that
$\tau \bigl[ \exp ( U) \bigr](\omega) \leq \epsilon$,
where
$$
U = \lambda \bigl[ r_2(\theta) - r_1(\theta) \bigr] 
- \frac{\lambda^2}{2N^2} \sum_{i=1}^N \bigl( 
\chi_{i + N} - \chi_i \bigr)^2 + \log(\epsilon).
$$
Thus $\tau \Bigl\{ \exp \bigl[ \tau'(U) \bigr] \Bigr\} \leq 
\tau \circ \tau' \bigl[ \exp (U) \bigr] = \tau \bigl[ \exp(U) \bigr] \leq \epsilon$,
from the convexity of the exponential function and the fact that
$\tau'$ is a (regular) conditional probability measure.

But \hfill
$\ds
\tau'(U) = \lambda \bigl[ r_2(\theta) - r_1(\theta) \bigr] 
- \frac{\lambda^2}{2 N} \tau' (V)
+ \log(\epsilon^{-1}),
$ \hfill {} \\
where $V = \frac{1}{N} \sum_{i=1}^N(\chi_{i+N} - \chi_i)^2$.
Noticing that
\begin{multline*}
\tau'(V) = \frac{1}{N} \sum_{i=1}^N (\chi_i + \chi_{i+N})
- 2  \left( \frac{1}{N}\sum_{i=1}^N \chi_i \right) \left( 
\frac{1}{N} \sum_{i=1}^N \chi_{i+N} \right) \\ = 
r_1(\theta) + r_2(\theta) - 2 r_1(\theta) r_2(\theta),
\end{multline*}
we get
\begin{lemma}
\label{lemma3.1}
For any exchangeable random variable $\eta$, 
\begin{multline*}
\tau \biggl\{ \exp \biggl[ \lambda \bigl[ r_2(\theta) - r_1(\theta) \bigr] 
- \frac{\lambda^2}{2N} \Bigl[ r_1(\theta) + r_2(\theta) - 2 r_1(\theta)
r_2(\theta) \Bigr] - \eta \biggr] \biggr\} (\omega)
\\ \leq \exp (- \eta) (\omega), \qquad \omega \in \Omega.
\end{multline*}
\end{lemma}

As a consequence,
\begin{lemma}
For any exchangeable probability distribution 
$\PP$, any exchangeable prior distribution $\pi$, with 
$\PP$ probability at least $1 - \epsilon$, for any $\theta \in \Theta$,
$$
r_2(\theta) \leq r_1(\theta) + \frac{\lambda}{2N}
\Bigl[ r_1(\theta) + r_2(\theta) - 2 r_1(\theta) r_2(\theta) \Bigr] + 
\frac{d}{\lambda},
$$
where $d = \inf \Bigl\{ - \log \bigl[ \pi(\theta') \epsilon \bigr]
: \theta' \in \Theta, \Psi(\theta') = \Psi(\theta) \Bigr\}$.
\end{lemma}
\begin{rmk}
As a special case, we can take $d = \log \bigl[ \C{N}(X_1^{2N}) \bigr]
- \log(\epsilon)$. This corresponds to the case when $\pi$ is chosen to
be the uniform distribution on $\Theta'$, using the remark that each 
$f_{\theta}$, $\theta \in \Theta$ coincides with some $f_{\theta'}$, 
$\theta' \in \Theta'$ on the design $\{X_i: i=1, \dots, 2N \}$.
\end{rmk}

We would like to prove a little more, showing that it is legitimate
to take in the previous equation
$$
\lambda = \left( \frac{2 N d 
}{r_1(\theta) + r_2(\theta) - 2 r_1(\theta) r_2(\theta)} \right)^{1/2} = 
\sqrt{\frac{2Nd}{\tau'(V)}}.
$$
This is not so clear, since this quantity is not (even partially)
exchangeable. Anyhow we can write the following:
\begin{multline*}
\sqrt{\frac{2Nd}{\tau'(V)}} \bigl\lvert r_2(\theta) - r_1(\theta) \bigr\rvert 
\leq \tau'(V^{-1/2}) \sqrt{2Nd} \bigl\lvert r_2(\theta) - r_1(\theta) \bigr\rvert 
\\ = \tau' \left( \sqrt{\frac{2Nd}{V}} \bigl\lvert r_2(\theta)-r_1(\theta)
\bigr\rvert \right),
\end{multline*}
because $r \mapsto r^{-1/2}$ is convex. Moreover, using successively
the fact that $\tau'(V)$ is a symmetric function of $r_1(\theta)$ and $r_2(\theta)$,
the fact that $\cosh$ is an even function, the previous inequality, 
the convexity of $\cosh$, the invariance $\tau = \tau \circ \tau'$, the 
invariance of $V$ under $\omega \mapsto \bigcirc_{i=1}^N \tau_i^1(\omega)$,
and the fact that $V$ is almost surely constant under each $\tau_i$,
we get the following chain of inequalities:
\pagebreak
\begin{multline*}
\tau \biggl\{ \exp \biggl[ \sqrt{\tfrac{2Nd}{\tau'(V)}}\bigl[ r_2(\theta) - 
r_1(\theta) \bigr] - d + \log(\epsilon) \biggr] \biggr\} 
\\ = 
\tau \biggl\{ \cosh \biggl[ \sqrt{\tfrac{2Nd}{\tau'(V)}}\bigl[ r_2(\theta) - 
r_1(\theta) \bigr] \biggr] \biggr\} \exp \bigl[- d + \log(\epsilon) \bigr]
\\ =
\tau \biggl\{ \cosh \biggl[ \sqrt{\tfrac{2Nd}{\tau'(V)}}\bigl\lvert r_2(\theta) - 
r_1(\theta) \bigr\rvert \biggr] \biggr\} \exp \bigl[- d + \log(\epsilon) \bigr]
\\ \leq \tau \biggl\{ \cosh \biggl[ \tau' \Bigl[ 
\sqrt{\tfrac{2Nd}{V}}\bigl\lvert r_2(\theta) - 
r_1(\theta) \bigr\rvert \Bigr] \biggr] \biggr\} \exp \bigl[- d + \log(\epsilon) \bigr]
\\ \leq \tau \biggl\{ \tau' \biggl[ \cosh \Bigl[ 
\sqrt{\tfrac{2Nd}{V}}\bigl\lvert r_2(\theta) - 
r_1(\theta) \bigr\rvert \Bigr] \biggr] \biggr\} \exp \bigl[- d + \log(\epsilon) \bigr]
\\ = \tau \biggl\{ \cosh \Bigl[ 
\sqrt{\tfrac{2Nd}{V}}\bigl[ r_2(\theta) - 
r_1(\theta) \bigr] \Bigr] \biggr\} \exp \bigl[- d + \log(\epsilon) \bigr]
\\ = \tau \biggl\{ \exp \biggl[ 
\sqrt{\tfrac{2Nd}{V}}\bigl[ r_2(\theta) - 
r_1(\theta) \bigr] - d + \log(\epsilon) \biggr] \biggr\}
\\ = \tau \circ \bigcirc_{i=1}^N \tau_i \biggl\{ \exp \biggl[ 
\sqrt{\tfrac{2Nd}{V}}\bigl[ r_2(\theta) - 
r_1(\theta) \bigr] - d + \log(\epsilon) \biggr] \biggr\}
\leq \epsilon.
\end{multline*}

Thus with $\PP$ probability at least $1 - \epsilon$, for any $\theta \in 
\Theta$, 
$$
r_2(\theta) \leq r_1(\theta) + \sqrt{\frac{2 d \tau'(V)}{N}}\\
= r_1(\theta) + \sqrt{\frac{2d \bigl[ r_1(\theta) + r_2(\theta)
-2 r_1(\theta) r_2(\theta)\bigr]}{N}}.
$$
Solving this inequality in $r_2(\theta)$ ends the proof of theorem \ref{th3.3}.
\end{proof}

In the general case when $\PP$ is only partially exchangeable and
$k$ is arbitrary, we will obtain the following
\begin{thm}
\label{th3.4}
Let $d = \inf \Bigl\{ - \log \bigl[ \pi(\theta) \epsilon \bigr]
: \theta' \in \Theta, \Psi(\theta') = \Psi(\theta) \Bigr\}$ and
\begin{multline*}
B(\theta) = \left( 1 + \frac{2d}{N}\right)^{-1}
\bigg\{ r_1(\theta) + \frac{d}{N} \Bigl\{ 1 + k^{-1}\bigl[
1 - 2 r_1(\theta) \bigr] \Bigr\} \\
+ (1 + k^{-1})\sqrt{\frac{2d}{N}r_1(\theta)\bigl[ 1 - r_1(\theta) 
\bigr] + \frac{d^2}{N^2}} \biggr\}.
\end{multline*}
For any partially exchangeable probability distribution 
$\PP$, with $\PP$ probability at least $1 - \epsilon$, for any $\theta 
\in \Theta$ such that $r_1(\theta) < 1/2$ and  $B(\theta) \leq 1/2$, 
$r_2(\theta) \leq B(\theta)$.
\end{thm}

As a special case, the theorem holds with
$d = \log \bigl[ \C{N}\bigl(X_1^{(k+1)N}\bigr) \bigr] + \log(\epsilon^{-1})$.
When using a set of binary classification rules 
$\{f_{\theta} : \theta \in \Theta\}$
whose VC dimension is not greater than $h$, we can use the bound 
$\ds d \leq h \log \left( \frac{e (k+1)N}{h} \right) - \log(\epsilon)$.
The result is satisfactory when $k$ is large, because in this case
$(1 + k^{-1})$ is close to one. This will be useful in the inductive
case.

Let us carry on our numerical example in the binary classification
case: taking $N = 1000$, $h = 10$,
$\epsilon = 0.01$ and 
$r_1(\theta) = 0.2$, we get a bound $B(\theta) \leq 0.4203$ for values
of $k$ ranging from $15$ to $18$, showing that increasing the size
of the shadow sample has an increased impact when the improved variance
term is used.

\begin{proof}
Let $\Phi(p) = (p \wedge \tfrac{1}{2}) \bigl[ 1 - (p \wedge \tfrac{1}{2})\bigr]$.
This is obviously a concave function. We have proved that
$$
\biggl( \bigcirc_{i=1}^N \tau_i \biggr) \biggl\{ 
\exp \Bigl[ \lambda \bigl[ r_2(\theta) - r_1(\theta) \bigr] - \eta
\Bigr] \biggr\} 
\leq \exp \biggl[ \frac{(1 + k^{-1})^2 \lambda^2}{2 N^2} \sum_{i=1}^N
\Phi(p_i) - \eta \biggr].
$$
As 
$$
\frac{1}{N} \sum_{i=1}^N \Phi(p_i) \leq \Phi \left( 
\frac{1}{N} \sum_{i=1}^N p_i \right) = \Phi \left( 
\frac{r_1(\theta) + k r_2(\theta)}{k+1} \right), 
$$
this shows that 
\begin{multline*}
\biggl( \bigcirc_{i=1}^N \tau_i \biggr) \biggl\{ 
\exp \Bigl[ \lambda \bigl[ r_2(\theta) - r_1(\theta) \bigr] - \eta
\Bigr] \biggr\} \\ \leq 
\exp \biggl[ \frac{ (1+k^{-1})^2 \lambda^2 }{2N} \Phi \biggl( 
\frac{r_1(\theta) + k r_2(\theta)}{k+1} \biggr) - \eta \biggr].
\end{multline*}

Taking $\ds \eta = \frac{(1+k^{-1})^2\lambda^2}{2N} \Phi \left( 
\frac{r_1(\theta)+ kr_2(\theta)}{k+1} \right) - \log(\epsilon)$,
and
$$
\lambda = \left( \frac{2Nd}{(1+k^{-1})^2 \Phi\left(
\frac{r_1(\theta)+kr_2(\theta)}{k+1}\right)}\right)^{1/2},
$$
where $d = \inf \bigl\{ - \log\bigl[\pi(\theta) \epsilon \bigr]
: \theta' \in \Theta, \Psi(\theta') = \Psi(\theta) \Bigr\}$,
we get that
with $\PP$ probability at least $1 - \epsilon$,
for any $\theta \in \Theta$,
$$
r_2(\theta) - r_1(\theta) \leq \left( \frac{2(1+k^{-1})^2 \Phi
\left( \frac{r_1(\theta)+kr_2(\theta)}{k+1} \right) d}{N} \right)^{1/2}.
$$

Solving this inequality in $r_2(\theta)$ ends the proof of theorem \ref{th3.4}.
\end{proof}

\section{The inductive setting}

We will integrate with respect to $\PP(\cdot|Z_1^N)$ theorem 
\ref{th2.1} and its variants. Let us start with theorem \ref{th3.4}.
Let us consider the non identically distributed independent case,
assuming thus that $\PP = \Bigl[ \bigotimes_{i=1}^N
P_i \Bigr]^{\otimes (k+1)}$.\\[10pt]
Let \hfill $\ds R(\theta) = \frac{1}{N} \sum_{i=1}^N P_i\bigl[Y_i 
\neq f_{\theta}(X_i)\bigr]$ \hfill {} \\[10pt] and \hfill 
$\ds \rr(\theta) = \frac{r_1(\theta) + k r_2(\theta)}{k+1}.$ \hfill {} \\[10pt]
Let \hfill $\PP'(h) = \PP(h \lvert Z_1^N).$ \hfill {} 
\newcommand{\piB}{\Bar{\pi}}

\pagebreak
\begin{lemma} For any partially exchangeable prior distribution $\pi$,
any partially exchangeable positive function 
$\zeta : \Theta \rightarrow \RR_+^*$,
\begin{multline*}
\PP \biggl\{ \sup_{\theta \in \Theta} \int_{\lambda = 0}^{+ \infty}
\zeta \exp \biggl[ \lambda \bigl[ R(\theta) - r_1(\theta) - \zeta 
\bigr] \\ - \frac{(1 + k^{-1})^2 \lambda^2}{2N}  \Phi \biggl( 
\frac{r_1(\theta) + k R(\theta)}{k+1} \biggr) + \PP' 
\Bigl[ \log \bigl[\piB(\theta) \epsilon \bigr] \Bigr] \biggr] 
d \lambda \biggr\}  \leq \epsilon,
\end{multline*}
where $\piB(\theta) = \sup \Bigl\{ \pi(\theta') : \theta' \in \Theta, 
\Psi(\theta') = \Psi(\theta) \Bigr\}$.
\end{lemma}
\begin{proof}
Let 
$$
U' = \lambda \bigl[ R(\theta) - r_1(\theta) - \zeta 
\bigr] \\ - \frac{(1 + k^{-1})^2 \lambda^2}{2N}  \Phi \biggl( 
\frac{r_1(\theta) + k R(\theta)}{k+1} \biggr) + \PP' 
\Bigl[ \log \bigl[\piB(\theta) \epsilon \bigr] \Bigr].
$$
Let 
$$
U = \lambda \bigl[ r_2(\theta) - r_1(\theta) - \zeta 
\bigr] \\ - \frac{(1 + k^{-1})^2 \lambda^2}{2N}  \Phi \bigl[ 
\rr(\theta)\bigr]+ \log \bigl[\piB(\theta) \epsilon \bigr].
$$
The function $\Phi$ being concave,
$$
\zeta \exp ( U' ) \leq \zeta \exp \bigl[ \PP' (U) \bigr]
\leq \PP' \bigl[ \zeta \exp(U) \bigr]. 
$$
Thus 
\begin{multline*}
\sup_{\theta} \int_{\lambda=0}^{+ \infty} \zeta \exp(U') d\lambda 
\leq \sup_{\theta \in \Theta} \int_{\lambda=0}^{+ \infty} 
\PP' \bigl[ \zeta \exp(U) \bigr] d \lambda
\\ \leq \sup_{\theta \in \Theta} \PP'
\left( \int_{\lambda=0}^{+\infty} \zeta \exp (U) d \lambda \right) 
\leq \PP' \left( \int_{\lambda = 0}^{+ \infty} \sup_{\theta 
\in \Theta} \bigl[ \zeta \exp (U) \bigr] d \lambda
\right)  
\end{multline*}
Moreover
$$
\sup_{\theta \in \Theta} \bigl[ \zeta \exp (U) \bigr] 
\leq \pi \bigl[ \zeta \exp (S) \bigr],
$$
where
$$
S = U - \log \bigl[ \pi(\theta) \bigr] = \lambda \bigl[ r_2(\theta) - r_1(\theta) - \zeta 
\bigr] \\ - \frac{(1 + k^{-1})^2 \lambda^2}{2N}  \Phi \bigl[ 
\rr(\theta)\bigr]+ \log (\epsilon ).
$$
Thus
\begin{multline*}
\PP \left( \sup_{\theta \in \Theta} \int_{\lambda = 0}^{+ \infty}
\zeta \exp ( U' ) d \lambda \right) 
\\ \leq \PP \left[ \PP' \left( \int_{\lambda=0}^{+ \infty} 
\pi \bigl[ \zeta \exp (S) \bigr] d \lambda \right) \right]\\
= \PP \left( \int_{\lambda =0}^{+ \infty} \pi \bigl[ 
\zeta \exp (S) \bigr] d \lambda \right) \\ = 
\PP \left[ \Bigl( \bigcirc_{i=1}^N \tau_i \Bigr)  
\left( \int_{\lambda = 0}^{+ \infty} \pi \bigl[ \zeta \exp (S) 
\bigr] d \lambda \right) \right] \\ 
= \PP \left\{ \pi \left[ \int_{\lambda = 0}^{+ \infty} 
\zeta \Bigl( \bigcirc_{i=1}^N \tau_i \Bigr)
\bigl[ \exp (S) \bigr] d \lambda \right] \right\}.
\end{multline*}
But we have established on the occasion of the proof of theorem
\ref{th3.4} that 
$$
\Bigl( \bigcirc_{i=1}^N \tau_i \Bigr)
\bigl[ \exp (S) \bigr] \leq \epsilon \exp(-\zeta \lambda).
$$
This proves that
$$
\PP \left( \sup_{\theta \in \Theta} \int_{\lambda = 0}^{+ \infty} 
\zeta \exp (U') d \lambda \right) 
\leq \epsilon,
$$
as stated in the lemma.
\end{proof}

\begin{thm}
\label{th4.2}
Let 
$$
B(\theta) = \left( 1 + \frac{2 d'}{N} \right)^{-1} 
\biggl\{ r_1(\theta) + \frac{d'}{N} 
+ \sqrt{ \frac{2 d' r_1(\theta) \bigl[ 1 - r_1(\theta) \bigr] }{N} 
+ \frac{{d'}^2}{N^2}} \biggr\},
$$
where 
$$
d' = d (1 + k^{-1})^2 \left( 1 - \frac{\log(\alpha)}{2 d} + \frac{\alpha}{
\sqrt{\pi d}} \right)^2,
$$
and 
$$
d = - \PP \Bigl\{ \log \bigl[ \epsilon \Bar{\pi}(\theta) \bigr] \lvert Z_1^N
\Bigr\}. 
$$
Let us notice that it covers the case when
$$
d = \PP \Bigl\{ \log \bigl[ \C{N}(X_1^{(k+1)N}) \bigr] \lvert Z_1^{N}
\Bigr\} + \log(\epsilon^{-1}).
$$
In this case, when $\lvert \C{Y} \rvert = 2$ and 
the set of classification rules has a VC dimension not greater than $h$, 
$$
d \leq h \log \left( \frac{e (k+1)N}{h} \right) + \log(\epsilon^{-1}).
$$
With $\PP$ probability at least $1 - \epsilon$, for any $\theta \in 
\Theta$, $R(\theta) \leq B(\theta)$ when $r_1(\theta) < 1/2$ and 
$B(\theta) \leq 1/2$.
\end{thm}

In the case when the model has a VC dimension not greater than $h$,
we can bound as  mentioned in the theorem the random variable
$d$ with the constant 
$$
d^{\star} = h \log \left( \frac{e(k+1)N}{h}
\right) + \log(\epsilon^{-1}).
$$
We can then optimize the choice of
$\alpha$ by taking $\alpha = \frac{1}{2} \sqrt{\frac{\pi}{d^{\star}}}$. 
This leads to
$$
d' \leq d^{\star} ( 1 + k^{-1})^2 \left[ 1 + \frac{1}{2 d^{\star}}
\log\left( 2 e \sqrt{\frac{d^{\star}}{\pi}}\right)\right]^2.
$$
We can also approximately optimize 
$$
(1 + k^{-1})^2 \log \left( \frac{eN(k+1)}{h} \right)
$$
by taking $k = 2 \log \left( \frac{ e N}{h} \right)$.

Let us resume our numerical example to illustrate theorem
\ref{th4.2}.
Assume that $N = 1000$, $h = 10$ and $\epsilon = 10^{-2}$.
For $r_1(\theta) = 0.2$, we get
$B(\theta) \leq 0.4257$ for $k = 19$. 
More generally, we get
$$
B(\theta) \leq 0.828 \left\{ r_1(\theta) + 0.105 + 
\sqrt{ 0.209 \bigl[1 - r_1(\theta) \bigr] r_1(\theta)
+ 0.011} \right\}. 
$$

For comparison, Vapnik's corollary \ref{cor1}
in the same situation gives a bound greater than $0.610$, and 
therefore not significant (since a random classification has 
a better expected error rate of $0.5$).

\begin{proof}
Let 
\begin{align*}
V & = (1 + k^{-1})^2 \Phi \left( \frac{r_1(\theta) + k R(\theta)}{k+1} 
\right),\\
d & = - \PP' \Bigl\{ \log \bigl[ \piB(\theta) \epsilon \bigr] \Bigr\},\\
\Delta & = R(\theta) - r_1(\theta) - \zeta.
\end{align*}
Let us remark that
$$
U' = - \frac{V}{2N} 
\left( \lambda - \frac{N \Delta}{V} \right)^2 + \frac{N \Delta^2}{2 V} 
- d.
$$
Thus 
$$
\int_{\lambda = 0}^{+ \infty} 
\zeta \exp (U') d \lambda \geq \B{1} (\Delta \geq 0) W, 
$$
where 
$$
W = 
\sqrt{\frac{\pi N}{2 V}} \zeta \exp \left( 
\frac{N \Delta^2}{2 V} - d \right).
$$
Thus, according to the previous lemma, 
$$
\PP \Bigl\{ \sup_{\theta \in \Theta} 
\bigl[ \B{1}( \Delta \geq 0) W \bigr] \Bigr\} \leq \epsilon.
$$
This proves that with $\PP$ probability at least 
$1 - \epsilon$, 
$$
\sup_{\theta \in \Theta} \bigl[ \B{1}(\Delta \geq 0) W \bigr] 
\leq 1.
$$
Translated into a logical statement this says that
with $\PP$ probability at least $1 - \epsilon$, 
either $\Delta < 0$, or $\log ( W) \leq 0$.

Let $V' = (1 + k^{-1}) \Phi\bigl[ R(\theta) \bigr]$.
Consider setting $\zeta = \alpha \sqrt{\frac{2 V'}{\pi N}}$,
where $\alpha$ is some positive real number.

We have proved that with $\PP$ probability at least $1 - \epsilon$,
$$
\frac{N \Delta^2}{2 V} \leq d - \log(\alpha) + \frac{1}{2}\log 
\left( \frac{V}{V'} \right),
$$
when $\Delta \geq 0$. But $\Phi$ is increasing and when $\Delta \geq 0$,
$R(\theta) \geq r_1(\theta)$, thus in this case $V' \geq V$, and 
we can weaken and simplify our statement to
$$
\frac{N \Delta^2}{2 V'} \leq d - \log(\alpha).
$$ 

Equivalently, with $\PP$ probability at least $1 - \epsilon$, 
for any $\theta \in \Theta$, 
$$
R(\theta) - r_1(\theta) \leq \sqrt{\frac{2 V' d}{N}}
\left( \sqrt{1 - \frac{\log(\alpha)}{d}} + \frac{\alpha}{\sqrt{\pi d}} 
\right)
$$
Using the fact that $\sqrt{1 + x} \leq 1 + \frac{x}{2}$,
we get that
$$
\bigl[ R(\theta) - r_1(\theta) \bigr]^2 \leq 
\frac{2 d' \Phi\bigl[R(\theta) \bigr]}{N},
$$
where $\ds d' = d (1 + k^{-1})^2 \left( 1 - \frac{\log(\alpha)}{2d} 
+ \frac{\alpha}{\sqrt{\pi d}} \right)^2$.
Since $\Phi(R) = R(1-R)$ when $R \leq 1/2$, this can be solved
in $R(\theta)$ in this case to end the proof of theorem \ref{th4.2}. 
\end{proof}

With a little more work we could have kept
$$
\frac{N \Delta^2}{2V} \leq d - \log(\alpha),
$$
leading to 
\begin{align*}
R(\theta) - r_1(\theta) & \leq \sqrt{\frac{2 V \bigl[d - 
\log(\alpha)\bigr]}{N}} + \alpha \sqrt{\frac{2 V'}{\pi N}},\\
\text{and }
\bigl[R(\theta) - r_1(\theta) \bigr]^2 
& \leq \frac{2V\bigl[ d - \log(\alpha)\bigr]}{N} 
+ 4 \alpha \frac{V'}{N} \sqrt{
\frac{d - \log{\alpha}}{\pi}} + \alpha^2 \frac{2V'}{\pi N}.
\end{align*}

This leads to the following
\begin{thm}
\label{th4.3}
Let us put 
\begin{align*}
c & = 2 \alpha \sqrt{\frac{d - \log(\alpha)}{\pi}} + \frac{\alpha^2}{\pi},\\
d_1 & = d - \log(\alpha) + (1 + k^{-1})^2 c,\\
d_2 & = (1 + k^{-1}) \Bigl\{ 
\bigl[d - \log(\alpha) \bigr] \bigl[ 1 - 
\tfrac{2 r_1(\theta)}{1 + k} \bigr] + (1 + k^{-1}) c \Bigr\},\\
d_3 & = (1 + k^{-1})^2 \Bigl\{ d - \log(\alpha) + c + 2 \tfrac{c}{N k^2}
\bigl[ d - \log(\alpha) \bigr] \Bigr\},\\
d_4 & = (1 + k^{-1}) \bigl[ d - \log(\alpha) + (1 + k^{-1}) c \bigr].
\end{align*}
Theorem \ref{th4.2} still holds when the bound $B(\theta)$ is strengthened to 
$$
B(\theta) = \left(1 + \frac{2 d_1}{N}\right)^{-1}
\left\{ r_1(\theta) + \frac{d_2}{N} + 
\sqrt{ \frac{2 d_3 r_1(\theta) \bigl[ 1 - r_1(\theta) \bigr]}{N} 
+ \frac{{d_4}^2}{N^2}} \right\}.
$$
\end{thm}

On the previous numerical example ($N = 1000$, $h = 10$, $\epsilon = 10^{-2}$,
$k = 19$, $\alpha = \frac{1}{2}\sqrt{\frac{\pi}{d^{\star}}}$, $r_1(\theta) = 
0.2$), we get a bound $B(\theta) \leq 0.4248$, instead of $B(\theta) \leq
0.4257$, showing that the improvement brought to theorem \ref{th4.2}
is not so strong, and therefore that theorem \ref{th4.2} is a satisfactory
approximation of theorem \ref{th4.3}.

Starting from lemma \ref{lemma2.1}, we can make the same kind of
computations taking $V = (1 + k^{-1}) \frac{r_1(\theta) + k R(\theta)}{k+1}$,
to obtain that with $\PP$ probability at least $1 - \epsilon$, 
$$
R(\theta) - r_1(\theta) \leq \sqrt{\frac{2 V' d}{N}} \left(1 
- \frac{\log(\alpha)}{2d} + \frac{\alpha}{\sqrt{\pi d}} \right), 
$$
where $V' = (1 + k^{-1})R(\theta)$. This proves the following
\begin{thm}
For any positive constant $\alpha$, with $\PP$ probability 
at least $1 - \epsilon$, for any $\theta \in \Theta$, 
$$
R(\theta) \leq r_1(\theta) + \frac{d'}{N} + \sqrt{\frac{2 d'r_1(\theta)}{N}
+ \frac{{d'}^2}{N^2}},
$$
where $d' = (1 + k^{-1}) \left(
1 - \frac{\log(\alpha)}{2d} + \frac{\alpha}{\sqrt{\pi d}}
\right)^2 d$.
\end{thm}

Our previous numerical application gives in this case a non significant bound
$R(\theta) \leq 0.516$, (for the best value of $k = 9$), showing that the improvement of the variance term
has a decisive impact when $r_1(\theta)$ is not small.

In the fully exchangeable case, when $k = 1$, a slightly better
result can be obtained, using lemma \ref{lemma3.1},
and thus putting 
\begin{align*}
V & = r_1(\theta) + R(\theta) - 2 r_1(\theta) R(\theta),\\
V' & = 2 R(\theta) \bigl[ 1 - R(\theta) \bigr].
\end{align*}
It leads to the following theorem

\begin{thm}
Let $d' = d - \log(\alpha)$ and
$$
c = 2 \alpha \sqrt{ \frac{ d - \log(\alpha)}{\pi}} + 
\frac{\alpha^2}{\pi}.
$$
Theorem \ref{th4.2} still holds when the bound is tightened 
to 
\begin{multline*}
B(\theta) = \left( 1 + \frac{4 c}{N}\right)^{-1} 
\Biggl\{ r_1(\theta) + \bigl[ 1 - 2 r_1(\theta) \bigr] 
\frac{d'}{N} + \frac{2 c}{N} 
\\ + \sqrt{ \frac{4 (d'+c)r_1(\theta) \bigl[1 - r_1(\theta) \bigr]}{N}
+ \frac{{d'}^2}{N^2} \bigl[ 1 - 2 r_1(\theta) \bigr]^2 
+ \frac{4 c (d'+c)}{N^2}} \Biggr\}.
\end{multline*}
\end{thm}
\begin{rmk}
Our previous numerical example gives in this case a bound
$B(\theta) \leq 0.460$, (for $\alpha = \frac{1}{2} \sqrt{\frac{\pi}{d}}$).
This shows that the improvement brought by a better variance term
is significant, but that the optimization of the size of the
shadow sample is also interesting.
\end{rmk}
\begin{rmk}
Note that we can take $\alpha = 1$. In this case, $d' = d$ and 
$c = 2 \sqrt{\frac{d}{\pi}} + \frac{1}{\pi}$. Note that we can
also take $\alpha = d^{-1/2}$, leading to
$d' = d + \frac{1}{2}\log(d)$ and 
$$
c = \frac{2}{\sqrt{\pi}} \sqrt{1 + \frac{\log(d)}{2d}} + 
\frac{1}{\pi d} \leq 
\frac{2}{\sqrt{\pi}} \left( 1 + \frac{\log(d)}{4d} \right) + \frac{1}{\pi d}
\leq \frac{3}{2}.
$$
\end{rmk}
\begin{rmk}
Note also that the bound can be weakened and simplified to
$$
B(\theta) \leq r_1(\theta) + \bigl[ 1 - 2 r_1(\theta) \bigr] 
\frac{d''}{N} + \sqrt{ \frac{4 d'' r_1(\theta) \bigl[ 1 - r_1(\theta) 
\bigr]}{N} + \frac{{d''}^2}{N^2} \bigl[ 1 - 
2 r_1(\theta) \bigr]^2},
$$
where $d'' = d' + 2c$. Taking $\alpha = d^{-1/2}$ gives
$d'' \leq d + \frac{1}{2}\log(d) + 3$.
\end{rmk}

Another technical possibility to get inductive bounds is to 
choose some near optimal value for $\lambda$, instead of averaging
over some exponential prior distribution on $\lambda$. 
\newcommand{\dB}{\Bar{d}}

This leads to the following theorem
\begin{thm}
\label{th4.6}
Let 
\begin{align*}
\dB & = \PP \Bigl\{ \log \bigl[ \piB(\theta)^{-1} \epsilon^{-1} \bigr] \Bigr\},\\
d & = \PP \Bigl\{ \log \bigl[ \piB(\theta)^{-1} \epsilon^{-1} \bigr] \lvert 
Z_1^N\Bigr\},\\
d' & = \frac{1}{4}(1+k^{-1})^2(\dB+d)(1 + \frac{d}{\dB}).
\end{align*}
Theorem \ref{th4.2} still holds when the bound is tightened
to 
$$
B(\theta) = \left( 1 + \frac{2 d'}{N} \right)^{-1}
\Biggl\{ r_1(\theta) + \frac{d'}{N} 
+ \sqrt{ \frac{2 d' r_1(\theta) \bigl[ 1 - r_1(\theta)\bigr]}{N}
+ \frac{{d'}^2}{N^2}} \Biggr\}.
$$
Moreover, putting $d^{\star} = \ess \sup_{\PP} \log \bigl[ \piB(\theta)^{-1}
\epsilon^{-1} \bigr]$, $d'$ can be replaced with
$d^{\star}(1 + k^{-1})^2$ in the previous bound. In the case of a
VC class of dimension $h$, $d^{\star}$ can be bounded by  
$h \log(\frac{eN}{h})$.
\end{thm}

Following our numerical example ($N=1000$, $h=10$, $r_1(\theta) = 0.2$),
we get an optimal value of $B(\theta) \leq 0.4213$ for $k$ ranging from 
$17$ to $19$. This shows that in this case going from the transductive
setting to the inductive one was done with an insignificant loss
of $0.001$. Although making use of a rather cumbersome flavor
of entropy term in the general case, theorem \ref{th4.6} 
provides the tightest bound in the case of a VC class.

\begin{proof}
Starting from
\begin{multline*}
\PP \biggl\{ \sup_{\theta \in \Theta}
\exp \biggl[ \lambda(\theta) \bigl[ R(\theta) - r_1(\theta) \bigr] 
\\ - \frac{\lambda(\theta)^2}{2N} (1+k^{-1})^2 \Phi \left( \frac{r_1(\theta)+kR(\theta)}{
1 + k}\right) + \PP'\Bigl\{ \log \bigl[ \piB(\theta) \epsilon \bigr] 
\Bigr\} \biggr] \biggr\} 
\leq \epsilon,
\end{multline*}
we can choose
$$
\lambda = \left[ \frac{- 2 N \PP \Bigl\{ \log \bigl[ \piB(\theta) \epsilon \bigr]
\Bigr\}}{
(1+k^{-1})^2\Phi\bigl[R(\theta)\bigr]} \right]^{1/2}.
$$
We get with $\PP$ probability at least $1 - \epsilon$,
$$
\lambda \bigl[ R(\theta)-r_1(\theta) \bigr] 
\leq \dB 
\frac{\Phi \left( \frac{r_1(\theta) + k R(\theta)}{1+k} \right)}{\Phi\bigl[ R(\theta) \bigr]} 
+ d.
$$
We can then remark that whenever $R(\theta) \geq r_1(\theta)$,
then 
$\frac{\Phi \left( \frac{r_1(\theta) + k R(\theta)}{1+k} \right)}{\Phi\bigl[ 
R(\theta) \bigr]} \leq 1$, 
to get 
$$
R(\theta) - r_1(\theta) \leq (1 + k^{-1}) \frac{\dB+d}{\sqrt{\dB}} 
\sqrt{ \frac{\Phi\bigl[ R(\theta) \bigr]}{2 N} }.
$$
Solving this inequality in $R(\theta)$ ends the proof of theorem \ref{th4.6}.
\end{proof}

In the same way, in the fully exchangeable case, starting
from
\begin{multline*}
\PP \biggl\{ \sup_{\theta \in \Theta}
\exp \biggl[ \lambda(\theta) \bigl[ 
R(\theta) - r_1(\theta) \bigr] \\ - \frac{\lambda(\theta)^2}{2N} 
\bigl[ R(\theta) + r_1(\theta) - 2 r_1(\theta) R(\theta) \bigr] 
+ d \biggr] \biggr\} \leq \epsilon,
\end{multline*}
we can take
$$
\lambda(\theta) = \sqrt{\frac{N\dB}{R(\theta) \bigl[ 1 - R(\theta) \bigr]}},
$$
to get
\begin{thm}
Let $\ds d' = \frac{1}{2}\dB \left( 1 + \frac{d}{\dB} \right)^2$, 
and assume that $\PP$ is fully exchangeable. 
Theorem \ref{th4.2} still holds when the 
bound is tightened to 
$$
B(\theta) = \left( 1 + \frac{2 d'}{N} \right)^{-1}
\Biggl\{ r_1(\theta) + \frac{d'}{N} 
+ \sqrt{ \frac{2 d' r_1(\theta) \bigl[ 1 - r_1(\theta)\bigr]}{N}
+ \frac{{d'}^2}{N^2}} \Biggr\}.
$$
Moreover, putting $d^{\star} = \ess \sup_{\PP}
\log \bigl[ \piB(\theta)^{-1} \epsilon^{-1} \bigr]$, 
$d'$ can be replaced with $2 d^{\star}$ in the 
previous bound. In the case of a VC class of dimension
$h$, $d^{\star}$ can be bounded by $h \log \left(\frac{eN}{h} \right)$.
\end{thm}

Our numerical example ($N=1000$, $h = 10$, $\epsilon = 0.0.1$
and $r_1(\theta) = 0.2$), gives a bound $B(\theta) \leq 0.445$.

\section{Using relative bounds}

Relative bounds were introduced in the PhD thesis of our student Jean-Yves 
Audibert \cite{Audibert}.
Here we will use them to sharpen Vapnik's bounds when 
$r_1(\theta)$ and $N$ are large (a flavor of how large they 
should be is given in the numerical application at the end of this
section). Audibert showed that chaining relative bounds can be
used to remove $\log(N)$ terms in Vapnik bounds. Here, we will
generalize relative bounds to increased shadow samples and 
will use only one step of the chaining method (lest we would
spoil the constants too much, the price to pay being a trailing
$\log\bigl[\log(N)\bigr]$ term which anyhow behaves like a constant
in practice).

Let us assume that $\PP$ is partially exchangeable. Let
$\theta$, $\theta' \in \Theta$, and let 
\begin{align*}
\chi_i & = \B{1}\bigl[Y_i \neq f_{\theta}(X_i) \bigr]
- \B{1} \bigl[ Y_i \neq f_{\theta'}(X_i) \bigr],\\
r_1'(\theta, \theta') & = r_1(\theta) - r_1(\theta')
= \frac{1}{N} \sum_{i=1}^N \chi_i,\\
r_2'(\theta, \theta') & = r_2(\theta) - r_2(\theta') 
= \frac{1}{kN} \sum_{i=N+1}^{(k+1)N} \chi_i.
\end{align*}

For any real number $x$, let $g(x) = x^{-2}\bigl[ \exp(x) - 1 - x \bigr]$.
As it is well known, $x \mapsto g(x) : \RR \rightarrow \RR$ is 
an increasing function. This is the key argument in the proof of
Bernstein's deviation inequality.

Let 
$$
\ell(\theta, \theta') = \frac{1}{(k+1)N} 
\sum_{i=1}^{(k+1)N} \B{1} \bigl[ f_{\theta}(X_i) \neq f_{\theta'}(X_i) \bigr].
$$

\begin{lemma}
For any partially exchangeable 
random variable $\lambda : \Omega \rightarrow \RR$,
\begin{multline*}
\left( \bigcirc_{i=1}^N \tau_i \right) \exp 
\biggl\{ \lambda \bigl[ r_2'(\theta, \theta') - r_1'(\theta, \theta') 
\bigr] \\ - g \bigl[ (1 + k^{-1}) \tfrac{2 \lambda}{N} \bigr] 
(1 + k^{-1})^2 \frac{\lambda^2}{N} \ell (\theta, \theta') 
+ \log(\epsilon) \biggr\} \leq \epsilon.
\end{multline*}
\end{lemma}
\pagebreak
\begin{proof}
For any partially exchangeable random variable $\eta$, 
\begin{multline*}
\log \biggl\{ \left( \bigcirc_{i=1}^N \tau_i \right)
\exp \biggl[ \lambda \bigl[ r_2'(\theta, \theta') 
- r_1'(\theta, \theta') \bigr] - \eta \biggr] \biggr\} \\ = - \eta +  
\sum_{i=1}^N \log \biggl\{ \tau_i 
\exp \biggl[ \frac{\lambda}{N} \biggl( 
\frac{1}{k} \sum_{j=1}^k \chi_{i + jN} - \chi_i \biggr) \biggr] \biggr\}
\\ = - \eta + \sum_{i=1}^N 
\log \biggl\{ 
\exp \biggl( 
\frac{\lambda}{kN} \sum_{j=0}^N \chi_{i + jN} \biggr) 
\tau_i \exp \biggl( - (1 + k^{-1}) \frac{\lambda}{N} \chi_i \biggr) 
\biggr\}\\ 
= -\eta + \frac{\lambda}{kN} \sum_{i=1}^{(k+1)N} 
\chi_i + \sum_{i=1}^N \log \biggl\{ 
\tau_i \biggl[ \exp \Bigl( - (1 + k^{-1})
\frac{\lambda}{N} \chi_i \Bigr) \biggr] \biggr\}.
\end{multline*}

Now we can apply Bernstein's inequality to
$$
\log \biggl\{ \tau_i \exp \biggl[ - (1 + k^{-1}) \frac{\lambda}{N} \chi_i 
\biggr] \biggr\},
$$

to show that 
\begin{multline*}
\log \biggl\{ \tau_i \exp \biggl[ 
- (1 + k^{-1}) \frac{\lambda}{N} \chi_i \biggr] \biggr\} 
\leq \frac{\lambda}{kN} \sum_{j=0}^k \chi_{i + j N} 
\\ + (1 + k^{-1})^2 \frac{\lambda^2}{N^2} 
\tau_i \Bigl[ ( \chi_i - p_i)^2 \Bigr] g \left[ 
\frac{2 \lambda}{N} (1 + k^{-1}) \right],
\end{multline*}
where we have put $\ds p_i = \tau_i (\chi_i ) = \frac{1}{k+1} \sum_{j=0}^k 
\chi_{i + jN}$. 
Anyhow, let us reproduce the proof of this statement here, for the 
sake of completeness. Let us put $\alpha = (1 + k^{-1}) \frac{\lambda}{N}$. 
\begin{multline*}
\log \Bigl\{ \tau_i \Bigl[ \exp ( - \alpha \chi_i ) \bigr] 
\Bigr\}  = - \alpha p_i \\ + 
\log \biggl\{ 1 + \tau_i \biggl[ 
\exp \Bigl[ - \alpha (\chi_i - p_i ) \bigr] 
- 1 - \alpha ( \chi_i - p_i) \biggr\} 
\\ \leq - \alpha p_i + 
\tau_i \Bigl[ \alpha^2 ( \chi_i - p_i)^2 g \bigl[ 
- \alpha ( \chi_i - p_i ) \bigr] \Bigr]
\\ \leq - \alpha p_i + g(2 \alpha) \tau_i \Bigl[ \alpha^2 (\chi_i - p_i)^2 \Bigr] .
\end{multline*}
We can now use the bound $\tau_i \Bigl[ ( \chi_i - p_i )^2 \Bigr] 
\leq \tau( \chi_i^2)$ and remark that \linebreak
$\chi_i^2 \leq \B{1} \bigl[ f_{\theta}(X_i) \neq f_{\theta'}(X_i) \bigr]$,
to get
\begin{multline*}
\log \biggl\{ \left( \bigcirc_{i=1}^N \tau_i \right) 
\exp \biggl[ \lambda \bigl[ r_2'(\theta, \theta') - r_1'(\theta, \theta') 
\bigr] - \eta \biggr] \biggr\} \\ \leq g(2 \alpha) \alpha^2 \sum_{i=1}^N \tau_i 
\Bigl\{ \B{1}\bigl[ f_{\theta}(X_i) \neq f_{\theta'}(X_i) \bigr] \Bigr\} 
- \eta \\ = 
\frac{\lambda^2}{N} (1 + k^{-1})^2 g \left( \frac{2 \lambda(1+k^{-1})}{N}\right)
\ell(\theta, \theta') - \eta.
\end{multline*}
We end the proof by choosing
$$
\eta = g\left( \frac{2 (1 + k^{-1}) \lambda}{N} \right)
(1 + k^{-1})^2 \frac{\lambda^2}{N} \ell(\theta, \theta') - \log(\epsilon).
$$
\end{proof}

We deduce easily from the previous lemma the following
\begin{proposition}
For any partially exchangeable prior distributions $\pi$, \linebreak $\pi' : 
\Omega \rightarrow \C{M}_+^1(\Theta)$, for any partially exchangeable
probability measure $\PP \in \C{M}_+^1(\Omega)$, with $\PP$
probability at least $1 - \epsilon$, for any $\theta$, $\theta' 
\in \Theta$, 
$$
r_2'(\theta, \theta') - 
r_1'(\theta, \theta') \leq g\left( \frac{2 (1 + k^{-1}) \lambda}{N} \right)
(1 + k^{-1})^2 \frac{\lambda}{N} \ell(\theta, \theta')
- \frac{1}{\lambda} \log \bigl[ \piB(\theta) \piB'(\theta') 
\epsilon \bigr],
$$
where $\piB(\theta) = \sup \{ \pi(\theta'') : \theta'' \in \Theta, 
\Psi(\theta'') = \Psi(\theta) \}$, and an analogous definition is used
for $\piB'$.
\end{proposition}

Let us now assume that we use a set of binary classification rules $\{
f_{\theta} : \theta \in \Theta \}$ with VC dimension not greater than
$h$.

Let us consider in the following the values
$$
\xi_j = \frac{ \left\lfloor (k+1)N \exp(-j) \right\rfloor}{(k+1)N} 
\simeq \exp(-j),
$$
where $\lfloor x \rfloor$ is the lower integer part of the real number $x$.
Let us define
$$
d_j = h \log \left[ \frac{2 e^2 (k+1) N}{h \xi_j} \right] + \log 
\bigl[ e \log(N) (h+1) \epsilon^{-1} \bigr].
$$

\begin{proposition}
\label{prop5.3}
With $\PP$ probability at least $1 - \epsilon$, 
for any $\theta \in \Theta$, any \linebreak $j \in \{1, \dots, \lfloor\log(N)\rfloor \}$, 
there is $\theta_j' \in \Theta_j'$ such that 
$$
r_2(\theta) - r_1(\theta) \leq r_2(\theta_j') - r_1(\theta_j') 
+ \left[ g \left( \sqrt{\frac{8 d_j}{\xi_j N}} \right) + \frac{1}{2} 
\right] \sqrt{\frac{2 (1 + k^{-1})^2 \xi_j d_j}{N}}.
$$
\end{proposition}
\begin{proof}
Let us recall a lemma due to David Haussler \cite{Haussler} : when the VC dimension of 
$\{ f_{\theta} : \theta \in \Theta \}$ is not greater than $h$, then, 
for any $\xi = \frac{m}{(k+1)N}$, we can find some $\xi$-covering net
$\Theta_{\xi}' \subset \Theta$ for the distance $\ell$
(which is a random exchangeable object), such that 
$$
\lvert \Theta_{\xi}' \rvert \leq e (h+1) \left( \frac{2 e}{\xi} \right)^h.
$$

Let us put on 
$$
\bigsqcup_{j, 1 \leq j \leq \log(N)} \Theta_{\xi_j}'
$$ 
the prior probability distribution defined by 
$$
\pi'(\theta'_j) = \bigl( \lfloor \log(N) \rfloor \lvert \Theta'_{\xi_j} \rvert \bigr)^{-1}
\geq \left[ \log(N) e (h+1) \left( \frac{2e}{\xi_j} \right)^h \right]^{-1}, 
\qquad \theta_j' \in \Theta_j'.
$$

We see that with $\PP$ probability at least $1 - \epsilon$, for any 
$\theta \in \Theta$, any $j$, $1 \leq j \leq \log(N)$ 
there is $\theta_j' \in \Theta_{\xi_j}'$
such that $\ell(\theta, \theta_j') \leq \xi_j$, and therefore such that 
\pagebreak
\begin{multline*}
r_2(\theta, \theta_j') - r_1(\theta, \theta_j') 
\leq g \left[(1+k^{-1})\frac{2 \lambda}{N}\right] (1+k^{-1})^2
\frac{\lambda}{N} \xi_j \\ + \frac{1}{\lambda} 
\log \bigl[ \piB(\theta)^{-1} \piB'(\theta_j')^{-1} \epsilon^{-1} \bigr],
\end{multline*}
where we can take $\piB(\theta) \geq \left( \frac{e(k+1)N}{h} \right)^{-h}$
and where $\pi'(\theta_j')$ has been defined earlier.
We can then choose
$$
\lambda( \theta, \theta_j') = \left[
\frac{2 N \log \Bigl\{ \bigl[ \piB(\theta) \piB'(\theta_j') \epsilon \bigr]^{-1}
\Bigr\}}{
(1 + k^{-1})^2 \xi_j} \right],
$$
to prove proposition \ref{prop5.3}.
\end{proof}

On the other hand, from theorem \ref{th3.4} applied to $\bigsqcup \Theta_{\xi_j}'$ 
and $\pi'$, we see that with $\PP$ probability at least $1 - \epsilon$, 
for any $j$, $1 \leq j \leq \log(N)$ and any $\theta_j' \in \Theta_j'$, 
putting
$$
d_j' = h \log \left( \frac{2e}{\xi_j} \right) + \log \bigl[ 
e \log(N) (h+1) \bigr] - \log(\epsilon),
$$
we have
$$
r_2(\theta_j') - r_1(\theta_j') \leq 
\sqrt{ \frac{2}{N} ( 1+ k^{-1})^2 d_j' \Phi \left( 
\frac{r_1(\theta_j') + k r_2(\theta_j')}{1 + k} \right)}.
$$
We can then remark that
when $\ell(\theta, \theta_j') \leq \xi_j$, 
$$
\frac{1}{k+1} \bigl[
r_1(\theta_j') + k r_2(\theta_j') \bigr]
\leq \frac{1}{k+1} \bigl[ r_1(\theta) + k r_2(\theta) 
\bigr] + \ell(\theta, \theta_j') 
\leq \frac{1}{k+1} \bigl[ r_1(\theta) + k r_2(\theta) 
\bigr] + \xi_j.
$$
We have proved the following
\begin{thm}
\label{th5.4}
With $\PP$ probability at least $1 - 2\epsilon$, 
\begin{multline*}
r_2(\theta) - r_1(\theta) \leq \inf_{j \in \NN^*, 1 \leq j \leq \log(N)}
\left[ g\left( \sqrt{\frac{8 d_j}{\xi_j N}} \right) + \frac{1}{2} \right]
\sqrt{ \frac{2(1 + k^{-1})^2 \xi_j d_j}{N}} \\ + 
\sqrt{ \frac{2}{N}(1 + k^{-1})^2 d_j' \Phi
\left( \frac{r_1(\theta) + k r_2(\theta)}{1 + k} + \xi_j \right)}.
\end{multline*}
\end{thm}
\begin{rmk}
To use this theorem, we have to solve equations of the type
$$
r_2 - r_1 \leq a + b 
\left[ \Phi \left(\frac{r_1 + kr_2}{1 + k} + \xi \right) \right]^{1/2}.
$$
Whenever $r_1$ and the bound are less than $1/2$, this is equivalent
to
$$
r_2 \leq \frac{B + \sqrt{B^2 - AC}}{A},
$$
where
\begin{align*}
A & = 1 + \left( \frac{kb}{1+k} \right)^2,\\
B & = r_1 + a + \frac{k b^2}{2 (1 + k)^2} \bigl[ 
(1 + k)(1 - 2 \xi) - 2 r_1 \bigr],\\
C & = (r_1 + a)^2 - \frac{b^2}{(1 + k)^2} \bigl[ 
(1 + k) \xi + r_1 \bigr] \bigl[ (1 + k)(1 - \xi) - r_1 \bigr].
\end{align*}
\end{rmk}

Let us make some numerical application. We should take $N$ pretty large, 
because the expected benefit of this last theorem is to improve on 
the $\log(N)$ term (the optimization in $\xi_j$ allows to kill the
$\log(N)$ term in $d_j$ and be left only with $\log\bigl[\log(N)
\bigr]$ terms). So let us take $N = 10^6$, $h = 10$, $r_1(\theta) = 0.2$
and $\epsilon = 0.005$. For these values, theorem \ref{th3.4} gives a
bound greater than $0.2075$ and less than $0.2076$ when 
$k$ ranges from $24$ to $46$. Here we
obtain a bound less than $0.2070$ for $k$ ranging from $24$ to $46$,
the optimal values for $(k,j)$ being $(257,7)$, giving a bound less than
$0.20672$. The bound is less than $0.2068$ for $k$ ranging from $42$
to $19470$, showing that we can really use big shadow samples with 
theorem \ref{th5.4} !